\documentclass[a4paper]{article}

\usepackage{graphicx}
\usepackage{amsmath}
\usepackage{amsfonts}
\usepackage{amssymb}
\newcommand{\be}{\begin{equation}}

\newcommand{\ee}{\end{equation}}

\newcommand{\bma}{\begin{pmatrix}}
\newcommand{\ema}{\end{pmatrix}}

\newcommand{\ga}{\gamma}

\newcommand{\Om}{\Omega}

\newcommand{\om}{\omega}

\newcommand{\la}{\lambda}

\newcommand{\dz}{\wedge}

\newcommand{\R}{{\bf R }}

\newcommand{\ba}{\begin{array}}

\newcommand{\ea}{\end{array}}

\newcommand{\beq}{\begin{eqnarray}}

\newcommand{\eeq}{\end{eqnarray}}

\newtheorem{lm}{Lemma}

\newtheorem{thee}{Theorem}

\newtheorem{proo}{Proposition}

\newtheorem{co}{Corollary}

\newtheorem{rem}{Remark}

\newtheorem{deff}{Definition}

\newcommand{\bd}{\begin{deff}}

\newcommand{\ed}{\end{deff}}

\newcommand{\bl}{\begin{lm}}

\newcommand{\el}{\end{lm}}

\newcommand{\bp}{\begin{proo}}

\newcommand{\ep}{\end{proo}}

\newcommand{\bt}{\begin{thee}}

\newcommand{\et}{\end{thee}}

\newcommand{\bc}{\begin{co}}

\newcommand{\ec}{\end{co}}

\newcommand{\brm}{\begin{rem}}

\newcommand{\erm}{\end{rem}}

\newcommand{\der}{{\rm d}}

\newcommand{\tj}{\theta^1}
\newcommand{\td}{\theta^2}
\newcommand{\ttr}{\theta^3}
\newcommand{\tc}{(\Om_2-\bar{\Om}_2)}
\newcommand{\ap}{(\Om_2+\bar{\Om}_2)}
\begin{document}
\thispagestyle{empty}
\title{{3-dimensional Cauchy-Riemann structures\\ and\\ 
$2^{\underline{\rm nd}}$ order ordinary differential equations} 
\vskip 1.truecm
\author{
Pawe\l~ Nurowski\\
Instytut Fizyki Teoretycznej\\
Uniwersytet Warszawski\\
ul. Hoza 69, Warszawa\\
Poland\\
nurowski@fuw.edu.pl\\
\\
\\
George A J Sparling\\
Department of Mathematics\\
University of Pittsburgh\\
Pittsburgh PA\\
USA\\
sparling@twistor.org}
}
\maketitle

\begin{abstract}

The equivalence problem for second order ODEs given 
modulo point transformations
is solved in full analogy with the equivalence problem of
nondegenerate 3-dimensional CR structures. This approach enables 
an analog of the Feffereman
metrics to be defined. The conformal class of these 
(split signature) metrics is well
defined by each point equivalence class of second order ODEs. Its
conformal curvature is interpreted in terms of the basic 
point invariants of the corresponding class of ODEs.
\end{abstract}

\newpage

\noindent

\rm
\section{Introduction}
This paper aims to explain the relations between two classical
geometries: the geometry associated with 2nd order ordinary
differential equations defined modulo point transformations of
variables and the geometry of the 3-dimensional Cauchy-Riemann
structures.\\

\noindent 
The geometry associated with 2nd order ordinary differential equations,
considered modulo contact transformations, is trivial - all 2nd
order ODEs are locally contact equivalent to the equation $y''=0$. If
one considers (more natural) point transformations - all
diffeomorphisms of 
the plane $(x,y)$ - then their action on the
space of all 2nd order ODEs has nontrivial orbits - there exist 2nd
order ODEs that are not (even locally) point equivalent. An example of
such point inequivalent equations is given by $y''=0$ and $y''=y^2$.\\

\noindent 
In general, the equation $$y''=Q(x,y,y')$$ with the total differential 
$$D=\partial_x+y'\partial_y+Q\partial_{y'},$$
can be characterized by a number of relative invariants whose
vanishing or not is a point invariant property of the equation. The
two of these invariants of lowest order are
$$w_1=D^2Q_{y'y'}-4DQ_{yy'}-DQ_{y'y'}Q_{y'}+
4Q_{y'}Q_{yy'}-3Q_{y'y'}Q_y+6Q_{yy},$$
and
$$w_2=Q_{y'y'y'y'}.$$
These were known to S. Lie \cite{Lie} and used by
M. A. Tresse \cite{Tre1,tres} in his systematic study of an equivalence 
problem for second order ODEs given modulo point
transformations. E. Cartan, in his celebrated paper \cite{Cartpc} on projective
connections, used the class of 2nd order ODEs for
which the invariant $w_2$ vanished as an example of a geometry that
naturally give rise to a Cartan normal projective
connection\footnote{Cartan's observation has recently been understood
from the twistorial point of view in \cite{Hitchin} and generalized in
\cite{Newman}.}.\\

\noindent
The study of the geometry of Cauchy-Riemann (CR) structures was initiated by
H. Poincare \cite{poin}, who looked for a higher dimensional
generalization of the well known fact that two real analytic arcs in
{\bf C} are locally biholomorphically equivalent. Using a heuristic
argument he showed that generic two real 3-dimensional hypersurfaces
$N_1$ and $N_2$ embedded in ${\bf C}^2$ are not, even locally,
biholomorphically equivalent. This led B. Segre \cite{Segre} to study 
the equivalence problem for real hypersurfaces of codimension
1 in ${\bf C}^2$, given modulo the biholomorphisms, a problem which 
was later solved in full generality by E. Cartan
\cite{CartCR}. Generalization of the problem to ${\bf C}^n$ with $n>2$
led to the theory of CR-structures which is a part of several complex
variables theory and lies on the borders between analysis, geometry
and studies of PDEs. In this theory a particular role is played by the
conformal Fefferman metrics \cite{Fef} which are {\it Lorentzian} metrics
naturally defined on a circle bundle over each CR manifold. These
metrics were defined by Ch. Fefferman in 1976 and, surprisingly, 
were unnoticed by E. Cartan in his pioneering paper \cite{CartCR}.\\ 

\noindent
The appearence of Lorentzian metrics in the CR-structure theory is not
an accident. In the lowest dimension ($n=2$) this is due to the well
known fact \cite{rt1,rt2,tafel,trautopt} that 3-dimensional
Cauchy-Riemann structures are in
one-to-one correspondence with congruences of null geodesics without
shear in 4-dimensional space-times. Many physically interesting
space-times, such as Minkowski, Schwarzschild, Kerr-Newman, Taub-NUT,
Hauser, plane gravitational waves and Robinson-Trautman, etc., admit 
congruences of such geodesics. The understanding of space-times
admitting congruences of
shear-free and null geodesics from the point of
view of the corresponding CR-geometry has been quite fruitful in the
process of solving Einstein vacuum equations 
\cite{lewnur,lnt1,lnt2,NPracmag,lnt3}. In these papers the construction of
the solutions of the Einstein equations in terms of the CR-functions
of the corresponding CR geometry is very much in the
spirit of the Penrose's twistor theory \cite{Pen1,Pen2}. More
importantly, from quite
another point of view, space-times admitting shear-free congruences of
null geodesics are the Lorentzian analogs of Hermitian
geometries in 4-dimensions. Since I. Robinson played the
crucial role of introducing the shear-free property to General
Relativity A. Trautman has called such manifolds Robinson
manifolds \cite{NurTraut,t1,t2}.\\

\noindent
Although it is not immediately self-evident, the geometries
associated with 2nd order ODEs and the 3-dimensional CR structures
are closely related. This fact was known to B. Segre who, in
this context, was quoted in Cartan's paper \cite{CartCR}. Strangely
enough, Cartan in \cite{CartCR} only mentioned that such relations
existed but did not spend much time explaining what they
were. In addition he does not appear to have used these
relations to simplify his approach to the CR structures. In this paper
we explain Segre's observation in detail and reconstruct and devolop 
the results of the two Cartan's papers 
\cite{Cartpc,CartCR} from this point of view. \\

\noindent
Section 2 consists of two parts. The first part contains a review of
the concept of a 3-dimensional CR-structure. This is
defined as a natural generalization of the notion of classes of real 
3-dimensional hypersurfaces embedded, modulo biholomorphisms, in ${\bf
  C}^2$. Our definition is much more in the spirit of Cartan's
treatment of such 
hypersurfaces
than in the spirit of the modern theory of CR manifolds. This point of
view will be adapted in the whole paper. The first part of Section
2 ends with the quotation of Cartan's theorem (Theorem 1) solving the
equivalence problem for 3-dimensional CR structures. In the second
part of Section 2 we give the modern description of this theorem in
terms of Cartan's {\bf SU}(2,1) connection. We also show how Cartan
might have used his theorem to associate with each nondegenerate
3-dimensional CR structure the Fefferman class of metrics. We analyze
the Fefferman metrics using Cartan's normal conformal connection
associated with them, and give a new proof,based on the use of
Baston-Mason conditions \cite{BasMas}, of the fact \cite{Lewand} that these
metrics are conformal to Einstein metrics only if the curvature of
their CR structure's {\bf SU}(2,1) connection vanishes. We close this 
section with a formula for the Bach tensor
for the Fefferman metrics, expressed in terms of the curvature of the
corresponding CR structure's connection.\\

\noindent
Section 3, the main Section of the paper, explains the analogy between
3-dimensional CR structures and 2nd order ODEs defined modulo point
transformations. The basic ingredients of this analogy are given just
before Definition 3, which states what it means for two second
order ODEs to be point equivalent to each other. Comparison between
Definitions 2 and 3 makes the analogy self evident. Using this analogy
we are able to formulate Theorem 2 which solves the equivalence
problem for 2nd order ODEs given modulo point transformations. By the
analogy this theorem is literally the same as Theorem 1. The only
difference is that now the symbols appearing in the Theorem have
different interpretations. This new interpretation implies that behind the
equivalence problem for 2nd order ODEs modulo point transformations is
a certain Cartan {\bf SL}(3,{\bf R}) connection. This fact was, of
course, known to Cartan \cite{Cartpc}, but we are not sure if 
Cartan would present it in the spirit of our paper even if he
had a time machine at his disposal (Cartan's ODE paper \cite{Cartpc} 
dates from 1924, whereas his CR paper \cite{CartCR} is from 1932). After
Theorem 2 we give a local representation of the point invariants of a
2nd order ODE obtaining, in particular, Lie's basic relative
invariants $w_1$ and $w_2$. We proceed, exploiting the analogy, to
define an ODE analog of Fefferman metrics, which now have {\it split} 
signature. The conformal class of split signature metrics which is 
naturally associated with each 2nd order ODE given modulo point
transformations turn out to encode all the point invariant information
about the underlying class of ODEs. In particular, all the Cartan
invariants of the point equivalent class of ODEs are derived from the
Weyl curvature of the corresponding Fefferman-like metric. These
metrics are characterized by Proposition 1 and the Remark following
it, and, some time ago, were considered by one of us (GAJS) within the
general framework discused in \cite{Spar}. Unlike the CR structures
case there are point equivalent classes of 2nd order ODEs which have
nonvanishing curvature of the Cartan {\bf SL}(3,{\bf R}) connection
and for which the Fefferman-like metrics are conformal to Einstein
metrics. Such metrics may only correspond to the ODEs 
for which the Lie relative invariants $w_1$ and $w_2$
satisfy $w_1w_2=0$ and all of them are presented in the Appendix.\\

\noindent
Like all the Cartan invariants, the Lie invariants $w_1$ and $w_2$ are 
interpreted in terms of the Fefferman-like metrics associated with the 
class of 2nd order ODEs that defines them. It turns out that the
Fefferman metrics associated with a point equivalence class of 2nd
order ODEs is always of the algebraic type $N\times N'$ in the
Cartan-Petrov-Penrose \cite{c,pen,pet} classification of real-valued 
4-dimensional metrics. This means, in particular, that both the self-dual
and the anti-self-dual parts of their Weyl tensor have only one
independent component. It turns out that the self-dual part of this
tensor is proportional to $w_1$ and the anti-self-dual part is
proportional to $w_2$. Thus, the vanishing of one of Lie's relative invariants
makes the associated Fefferman-like metric half-flat. This partially
explains why, in such cases, these metrics may be conformal to
Einstein metrics. The rest of Section 3 is devoted to understanding
the fact that it is easy to find all 2nd
order ODEs for which $w_2=0$, since all of them are of the form
$$
y''=A_0(x,y)+A_1(x,y)y'+A_2(x,y)(y')^2+A_3(x,y)(y')^3,
$$ 
and it is quite hard to find $Q=Q(x,y,y')$ for which $w_1=0$. From the
Fefferman-like metrics point of view the switch bewteen $w_2$ and
$w_1$ is the switch between the self-dual and the anti-self-dual part
of their Weyl tensor. This suggest that invariants $w_1$ and $w_2$
should be on an equal footing of complexity. To see that this is indeed the
case requires another notion of duality - the duality between
the point equivalent classes of 2nd orderd ODEs. This duality was
mentioned by Cartan in \cite{Cartpc}. We explain it in detail at the
end of Section 3. In particular, in Proposition 2, and in the Example
preceeding it, we show how to construct solutions $Q=Q(x,y,y')$ of
$w_1=0$ knowing $Q$s which satisfy $w_2=0$. The understanding of this
duality in terms of the natural double fibration of the first jet
bundle associated with the ODE is also given.\\

\noindent
Finally, in Section 4 we give two applications of the theory presented
in Sections 2 and 3. The first consists of an algorithm for associating
a point equivalence class of 2nd order ODEs with a given 3-dimensional
CR structure. This may be of some use in General Relativity theory and 
may provide a new understanding of well known congruences of
shear-free and null geodesics. The second application is, as  far as
we know, the first example of a large class of split signature
4-metrics which satisfy the
Bach equations, are genuinely of algebraic type $N\times N'$ and are
not conformal to Einstein metrics.\\

\noindent
{\it Note on the conventions and the notation}\\

\noindent
We emphasize that in this paper all our considerations are purely
local and concerned with nonsingular points of the introduced
constructions. We also mention that, following the old tradition in
PDEs, we denote the partial derivatives with respect to the variable
associated with index $i$ by the corresponding subscript, e.g. 
$\frac{\partial G}{\partial z_i}=G_i$.

\section{3-dimensional CR structures}
A 3-dimensional CR structure is a structure which a 3-dimensional
hypersurface $N$ embedded in ${\bf C}^2$ acquires from the ambient
complex space. Following Elie Cartan \cite{CartCR} this structure can
be described in the language of differential forms as follows.\\

\noindent
Consider a 3-dimensional hypersurface $N$ in ${\bf C}^2$ defined by
means of a real function $G=G(z_1,z_2,\bar{z}_1,\bar{z}_2)$, such that
$G_1\neq 0$, via
$$
N=\{~(z_1,z_2)\in{\bf C}^2~|~G(z_1,z_2,\bar{z}_1,\bar{z}_2)=0~\}.
$$
All information about the structure acquired by $N$ from ${\bf C}^2$ 
can be encoded in the two 1-forms 
\be
\lambda=i(G_1\der z_1+G_2\der z_2)~~~~~{\rm and}~~~~\mu=\der z_2.\label{formy}
\ee
These forms have the following properties:
\begin{itemize}
\item $\lambda$ is real, $\mu$ is complex
\item $\lambda\dz\mu\dz\bar{\mu}\neq 0$ on $N$.
\end{itemize}
\noindent
Moreover, if $N$ underlies the biholomorphism 
$$z_1=z_1(z'_1,z'_2),~~~
  z_2=z_2(z'_1,z'_2)$$ 
the forms transform according to
$$
\lambda\to\lambda'=a\lambda~~~~{\rm and}~~~
\mu\to\mu'=b\mu+c\lambda,
$$ 
where $a\neq 0$ (real) and $b\neq 0, c$ (complex) are appropriate 
functions on $N$. It is easy to see that the vanishing of the 3-form 
$\lambda\dz\der\lambda$ is an invariant property under the
  biholomorphisms of ${\bf C}^2$. Thus, the two hypersurfaces 
$$N_1=\{~(z_1,z_2)\in {\bf C}^2 ~:~ z_1-\bar{z}_1=0~\}$$
and 
\be
N_2=\{~(z_1,z_2)\in {\bf C}^2 ~:~ |z_1|^2+|z_2|^2-1=0~\},\label{sfera}
\ee
with the corresponding forms $\lambda_1=i\der z_1$ and
$\lambda_2=i(\bar{z}_1\der z_1+\bar{z}_2\der z_2)$ are not
biholomorphically equivalent.\\

\noindent
The above considerations motivate an introduction of the following
structure on 3-manifolds.
\bd 
A CR-structure $[(\lambda,\mu)]$ on a 3-dimensional manifold $N$
is an equivalence class of pairs of 1-forms $(\lambda,\mu)$ such that
\begin{itemize}
\item $\lambda$ is real, $\mu$ is complex,
\item $\lambda\dz\mu\dz\bar{\mu}\neq 0$ on $N$
\item two pairs $(\lambda,\mu)$ and $(\lambda',\mu')$ are in the
  equivalence relation
  iff there exist functions $a\neq 0$ (real), $b\neq 0, c$ (complex)
  on $N$ such that 
$$
\lambda'=a\lambda,~~~~~~~
\mu'=b\mu+c\lambda,~~~~~~~
\bar{\mu}'=\bar{b}\bar{\mu}+\bar{c}\lambda.
$$
\end{itemize}
\noindent
A CR-structure is called nondegenerate iff
$$
\der\lambda\dz\lambda\neq 0;
$$
otherwise a CR-structure is degenerate.
\ed
An obvious class of examples of CR-structures is given by
biholomorphically equivalent classes of hypersurfaces in ${\bf
  C}^2$. The problem of classifying biholomorphically nonequivalent
hypersurfaces in ${\bf C}^2$ is therefore a part of the equivalence
problem of CR-structures.
\bd
Let $(N,[(\lambda,\mu)])$ and $(N',[(\lambda',\mu')])$ be 
two CR-structures on two 3-dimensional manifolds $N$ and $N'$. 
We say that $(N,[(\lambda,\mu)])$ and $(N',[(\lambda',\mu')])$ are
(locally) equivalent
iff, for any two representatives $(\lambda,\mu)\in [(\lambda,\mu)]$ and 
$(\lambda',\mu')\in[(\lambda',\mu')]$, there exists a (local) 
diffeomorphism $\phi: N\to N'$ and functions $a\neq 0$ (real), 
$b\neq 0,c$ (complex) on $N$ such that 
$$
\phi^*(\lambda')=a\lambda,~~~~~~~\phi^*(\mu')=b\mu+c\lambda,~~~~~~~\phi^*(\bar{\mu}')=\bar{b}\bar{\mu}+\bar{c}\lambda.
$$
\ed 
It is easy to see that all 3-dimensional {\it degenerate}
CR-structures are locally equivalent to the structure associated with
a biholomorphic class of hypersurfaces equivalent to the hypersurface
$N_1={\bf C}\times{\bf R}$. 
The equivalence problem for {\it nondegenerate} 3-dimensional CR-structures was solved
by Elie Cartan \cite{CartCR}. Given a nondegenerate CR-structure $(N,[(\lambda,\mu)])$ he considered the forms 
\be
\theta^1=b\mu+c\lambda,~~~~~~~\theta^2=\bar{b}\bar{\mu}+\bar{c}\lambda,~~~~~~~\theta^3=a\lambda,\label{theta}
\ee  
with some unspecified functions $a\neq 0$ (real), $b\neq 0$ and $c$
(complex). He viewed the forms as being well defined on an
8-dimensional space $P_0$ parametrized by the points of $N$ and by the
coordinates $(a,b,\bar{b},c,\bar{c})$. Using his equivalence method (see
e.g. \cite{Jac,Olver}) he then
constructed another 8-dimensional manifold $P$ on which the coframe
consisting of the forms $(\theta^1,\theta^2,\theta^3)$ and the five additional
well defined 1-forms $(\Om_2,\bar{\Om}_2,\Om_3,\bar{\Om}_3,\Om_4)$ constitued the system of
basic biholomorphic invariants of the CR-structure. More precisely, he
proved the following theorem.
\noindent
\bt
Every nondegenerate CR-structure $(N,[(\lambda,\mu)])$ 
uniquely defines an 8-dimensional manifold $P$, 1-forms
$\theta^1,\theta^2,\theta^3,\Om_2,\bar{\Om}_2,\Om_3,\bar{\Om}_3,\Om_4$
and functions  ${\cal R},\bar{\cal R},{\cal S},\bar{\cal S}$ on $P$
such that
\begin{itemize}
\item $\theta^1,\theta^2,\theta^3$ are as in (\ref{theta}),
  $\bar{\Om}_2,\bar{\Om}_3$ are respective complex conjugates of
  $\Om_2,\Om_3$ and $\Om_4$ is real,
\item
  $\theta^1\dz\theta^2\dz\theta^3\dz\Om_2\dz\bar{\Om}_2\dz\Om_3\dz\bar{\Om}_3\dz\Om_4\neq
  0$ at each point of $P$.
\end{itemize}
\noindent
The forms satisfy the following equations
\beq
&\der\theta^1=\Om_2\dz\theta^1+\Om_3\dz\theta^3\nonumber\\
&\der\theta^2=\bar{\Om}_2\dz\theta^2+\bar{\Om}_3\dz\theta^3\nonumber\\
&\der\theta^3=i\theta^1\dz\theta^2+(\Om_2+\bar{\Om}_2)\dz\theta^3\nonumber\\
&\der\Om_2=2
i\theta^1\dz\bar{\Om}_3+i\theta^2\dz\Om_3+\Om_4\dz\theta^3\nonumber\\
&\der\bar{\Om}_2=-2
i\theta^2\dz\Om_3-i\theta^1\dz\bar{\Om}_3+\Om_4\dz\theta^3\label{sys}\\
&\der\Om_3=\Om_4\dz\theta^1+\Om_3\dz\bar{\Om}_2+{\cal
  R}\theta^2\dz\theta^3\nonumber\\
&\der\bar{\Om}_3=\Om_4\dz\theta^2+\bar{\Om}_3\dz\Om_2+\bar{\cal
  R}\theta^1\dz\theta^3\nonumber\\
&\der\Om_4=i\Om_3\dz\bar{\Om}_3+\Om_4\dz(\Om_2+\bar{\Om}_2)+\bar{\cal
  S}\theta^1\dz\theta^3+{\cal S}\theta^2\dz\theta^3\nonumber.
\eeq
The functions ${\cal R}$, ${\cal S}$, and their respective 
complex conjugates $\bar{\cal R}$, $\bar{\cal S}$, satisfy  
\beq
&\der{\cal R}=-{\cal R}(\Om_2+3\bar{\Om}_2)-{\cal S}\theta^1+{\cal R}_2\theta^2+{\cal
  R}_3\theta^3\nonumber\\
&\der\bar{\cal R}=-\bar{\cal R}(\bar{\Om}_2+3\Om_2)-\bar{\cal S}\theta^2+\bar{\cal R}_2\theta^1+\bar{\cal
  R}_3\theta^3\label{funk}
\eeq
and
\beq
&\der{\cal S}=-{\cal S}(2\Om_2+3\bar{\Om}_2)-i{\cal
  R}\bar{\Om}_3+{\cal S}_1\theta^1+{\cal S}_2\theta^2+{\cal S}_3\theta^3\nonumber\\
&\der\bar{\cal S}=-\bar{\cal S}(2\bar{\Om}_2+3\Om_2)+i\bar{\cal
  R}\Om_3+\bar{\cal S}_2\theta^1+{\cal S}_1\theta^2+\bar{\cal S}_3\theta^3,\label{ghj}
\eeq
with appropriate functions ${\cal R}_2,{\cal R}_3,{\cal S}_1,{\cal
  S}_2,{\cal S}_3$ and their conjugates.\\
The function ${\cal S}_1$ satisfies 
\be 
{\cal S}_1=\bar{\cal S}_1.\label{bach}
\ee
\et

\noindent
The above theorem, stated in the modern language, means the
following. The manifold $P$ is a Cartan bundle $H\to P\to N$, with $H$
a 5-dimensional parabolic subgroup of ${\bf SU}(2,1)$. This latter group
preserves the $(2,1)$-signature hermitian form
$$
h(X,X)=
\bma
X^1,&X^2,&X^3
\ema
\hat{h}
\bma
\bar{X}^1\\
\bar{X}^2\\
\bar{X}^3
\ema, ~~~~~~ 
\hat{h}=
\bma 
0&0&2i\\
0&1&0\\
-2i&0&0
\ema.
$$ 
The forms
$\theta^1,\theta^2,\theta^3,\Om_2,\bar{\Om}_2,\Om_3,\bar{\Om}_3,\Om_4$
of the theorem can be collected into a matrix of 1-forms 
$$
\omega=
\begin{pmatrix}
\frac{1}{3}(2\Om_2+\bar{\Om}_2)&i\bar{\Om}_3&-\frac{1}{2}\Om_4\\
&&\\
\theta^1&\frac{1}{3}(\bar{\Om}_2-\Om_2)&-\frac{1}{2}\Om_3\\
&&\\
2\theta^3&2i\theta^2&-\frac{1}{3}(2\bar{\Om}_2+\Om_2),
\end{pmatrix}
$$
satisfying
$$
\om\hat{h}+\hat{h}\om^\dagger=0,
$$
which is an ${\bf su}(2,1)$-valued Cartan connection \cite{Kobayashi}
on $P$. It follows from
equations (\ref{sys}) that the curvature of
this connections is
$$
\Om=\der \om+\om\dz\om=
\begin{pmatrix}
0&i\bar{\cal R}\theta^1\dz\theta^3&-\frac{1}{2}\bar{\cal
  S}\theta^1\dz\theta^3-\frac{1}{2}{\cal S}\theta^2\dz\theta^3\\
&&\\
0&0&-\frac{1}{2}{\cal R}\theta^2\dz\theta^3\\
&&\\
0&0&0
\end{pmatrix}.
$$
It yields all the invariant information about the corresponding
CR-structure, very much in the way as the Riemann curvature yields all
the information about the Riemannian structure.\\

\noindent
{\bf Remark} Note that the assumption that ${\cal R}$ or ${\cal S}$
(and, therefore $\bar{\cal R}$ or $\bar{\cal S}$) are constant on $P$ 
is compatible with (\ref{funk}) iff ${\cal R}={\cal S}=0$ (and, 
therefore $\bar{\cal R}=\bar{\cal S}=0$). In such case the
curvature $\Om$ of the Cartan connection $\om$ 
vanishes, and it follows that
there is only one, modulo local equivalence, CR-structure with this property. It
coincides with the CR-structure, which the hypersurface $N_2={\bf
  S}^3$ acquires from the ambient space ${\bf C}^2$ via equations (\ref{formy}),
(\ref{sfera}).\\

\noindent
Using the matrix elements $\om^i_{~j}$ of the Cartan connection $\om$
it is convenient to consider the bilinear form
$$
G=-i\om^3_{~j}\om^j_{~1}.
$$
This form, when writen explicitly in terms of
$\theta^1,\theta^2,\theta^3,\Om_2,\bar{\Om}_2,\Om_3,\bar{\Om}_3,\Om_4$,
is given by
$$
G=2\theta^1\theta^2+\frac{2}{3i}\theta^3(\Om_2-\bar{\Om}_2).
$$
Introducing the basis of vector fields
$X_1,X_2,X_3,Y_2,\bar{Y}_2,Y_3,\bar{Y}_3,Y_4$, the respective dual of
$\theta^1,\theta^2,\theta^3,$ $\Om_2,\bar{\Om}_2,\Om_3,\bar{\Om}_3,\Om_4$,
one sees that $G$ is a form of signature $(+++-0000)$ with four 
degenerate directions corresponding to four vector fields $Z_I=(Y_2+\bar{Y}_2,Y_3,\bar{Y}_3,Y_4)$. These 
four directions span a
4-dimensional distribution which is integrable due to equations
(\ref{sys}). Thus, the Cartan bundle $P$ is foliated by 4-dimensional
leaves tangent to the degenerate directions of $G$. Moreover, 
equations (\ref{sys}) guarantee that 
$$
{\cal L}_{Z_I}~G=A_I~G,
$$ 
with certain functions $A_I$ on $P$,
so that the bilinear form $G$ is preserved up to a scale when Lie
transported along the leaves of the foliation. Therefore the
4-dimensional space $P/\hspace{-0.15cm}\sim$ of leaves of the foliation is naturally
equipped with a conformal class of Lorentzian metrics $[g_F]$, the
class to which the bilinear form $G$ naturally descends. The
Lorentzian metrics
\be
g_F=2\theta^1\theta^2+\frac{2}{3i}\theta^3(\Om_2-\bar{\Om}_2)\label{fef}
\ee
on $P/\hspace{-0.15cm}\sim$ coincide with the so called Fefferman 
metrics \cite{Fef} (see also \cite{Grah}) which Charles Fefferman associated with any nondegenerate 
CR-structure $(N,[(\lambda,\mu)]$.\\

\noindent 
Introducing the volume form 
$$
\eta=\frac{1}{3}\theta^1\dz\theta^2\dz\theta^3\dz(\Om_2-\bar{\Om}_2)
$$
on $P/\hspace{-0.15cm}\sim$ we observe that the Hodge dualization $*$
of the forms $\theta^3\dz\theta^1$ and $\theta^3\dz\theta^2$ read
$$
*(\theta^3\dz\theta^1)=-\frac{1}{i}(\theta^3\dz\theta^1)~~~~{\rm
  and}~~~~
*(\theta^3\dz\theta^2)=\frac{1}{i}(\theta^3\dz\theta^2).
$$
Thus $\theta^3\dz\theta^1$ is self-dual and $\theta^3\dz\theta^2$ is
anti-self-dual.\\

\noindent
A convenient way of analyzing the Fefferman metrics is to look for the
Cartan normal conformal connection associated with them. Given a nondegenerate
CR-structure $(N,[(\lambda,\mu)]$ we define an ${\bf so}(4,2)$-valued 
matrix of 1-forms $\tilde{\om}$ on $P$ via 
\be
\tilde{\om}=
\bma
\frac{1}{2}\ap&\frac{i}{2}\bar{\Om}_3&-\frac{i}{2}\Om_3&-\Om_4&\frac{i}{12}\tc&0\\
&&&&&\\
\tj&-\frac{1}{3}\tc&0&-\Om_3&\frac{i}{2}\tj&-\frac{i}{2}\Om_3\\
&&&&&\\
\td&0&\frac{1}{3}\tc&-\bar{\Om}_3&-\frac{i}{2}\td&\frac{i}{2}\bar{\Om}_3\\
&&&&&\\
\ttr&\frac{i}{2}\td&-\frac{i}{2}\tj&-\frac{1}{2}\ap&0&\frac{i}{12}\tc\\
&&&&&\\
\frac{1}{3i}\tc&\bar{\Om}_3&\Om_3&0&\frac{1}{2}\ap&-\Om_4\\
&&&&&\\
0&\td&\tj&\frac{1}{3i}\tc&\ttr&-\frac{1}{2}\ap
\ema.\label{18}
\ee
This is a pullback of the Cartan normal conformal connection
associated with the Fefferman
metric from the Cartan ${\bf SO}(4,2)$ conformal bundle to $P$. With
a slight abuse of the language we call $\tilde{\om}$ the Cartan
conformal coonnection. The
pulback of the curvature of this connection   
$$
\tilde{\Om}=\der\tilde{\om}+\tilde{\om}\dz\tilde{\om},
$$ 
is given by 
\beq
&\tilde{\Om}=\tilde{\Om}^++\tilde{\Om}^-=\nonumber\\
&\\
&=\bma
0&-\frac{i}{2}\bar{\cal R}&0&\bar{\cal S}&0&0\\
&&&&&\\
0&0&0&0&0&0\\
&&&&&\\
0&0&0&\bar{\cal R}&0&-\frac{i}{2}\bar{\cal R}\\
&&&&&\\
0&0&0&0&0&0\\
&&&&&\\
0&-\bar{\cal R}&0&0&0&\bar{\cal S}\\
&&&&&\\
0&0&0&0&0&0
\ema\theta^3\dz\theta^1+\bma
0&0&\frac{i}{2}{\cal
  R}&{\cal
  S}&0&0\\
&&&&&\\
0&0&0&{\cal R}&0&\frac{i}{2}{\cal R}\\
&&&&&\\
0&0&0&0&0&0\\
&&&&&\\
0&0&0&0&0&0\\
&&&&&\\
0&0&-{\cal
  R}&0&0&{\cal
  S}\\
&&&&&\\
0&0&0&0&0&0
\ema\theta^3\dz\theta^2.\nonumber
\eeq
Here $\tilde{\Om}^+$ and $\tilde{\Om}^-$ denote the self-dual 
and the anti-selfdual parts of $\tilde{\Om}$, respectively.\\

\noindent
The theory of the conformal connections \cite{Frit,Kobayashi,Koz,PN1} then implies that 
that the Weyl curvature 2-form $C$ of $g_F$ is given by\footnote{Here $C$
  has tensor indices $C^\mu_{~\nu}$, $\mu,\nu=1,2,3,4$ which are
  associated with the null tetrad
  $\theta^1,\theta^2,\theta^3,\theta^4=\frac{1}{3i}(\Om_2-\bar{\Om}_2)$ of $g_F$. In this tetrad $g_F=2\theta^1\theta^2+2\theta^3\theta^4$.}
\beq
&C=C^+ +C^-=\nonumber\\
&\label{20}\\
&=\bar{\cal R}\bma
0&0&0&0\\
&&&\\
0&0&1&0\\
&&&\\
0&0&0&0\\
&&&\\
-1&0&0&0
\ema\theta^3\dz\theta^1+{\cal R}\bma
0&0&1&0\\
&&&\\
0&0&0&0\\
&&&\\
0&0&0&0\\
&&&\\
0&-1&0&0
\ema\theta^3\dz\theta^2,\nonumber
\eeq
i.e. denoting the matrix elements of $\tilde{\Om}$ by $\tilde{\Om}^A_{~B}$,
$A,B=0,1,...5$, it is given by $\tilde{\Om}^A_{~B}$ with 
$A,B=1,2,3,4$. The very simple form of the Weyl curvature $C$ shows
that the Fefferman metric $g_F$ of any nondegenerate CR-structure
$(N,[(\la,\mu)])$ is of the Petrov type $N$.\\

\noindent
{\bf Remark} Note that the curvature of the Cartan
normal conformal connection $\tilde{\om}$ of the Fefferman metric $g_F$ yields
essentially the same information as the curvature of the 
${\bf su}(2,1)$-valued connection $\om$. This is due to the fact
\cite{BDS} that
$\om$ is simply an ${\bf su}(2,1)$ reduction of the Cartan normal
conformal connection associated with the Fefferman metric $g_F$.
In addition, this indicates the well
known fact that the Fefferman conformal class of metrics $[g_F]$ associated with
a given nondegenrate CR-structure $(N,[(\la,\mu)])$ yields all the
invariant information about $(N,[(\la,\mu)])$. In
particular, the ${\bf su}(2,1)$-curvature properties of the
CR-structure are totally encoded in the Weyl tensor 2-forms $C$ of the
corresponding Fefferman metrics. Note that although $C$
  explicitely involves only ${\cal R}$ and $\bar{\cal R}$, the  ${\cal
    S}$ and $\bar{\cal S}$ functions can be derived from them by means
  of equations (\ref{funk}). \\

\noindent
It is known \cite{Lewand} that the Fefferman metrics are conformal to
the Einstein metrics only in the case when the corresponding
CR-structure is flat ($\Om=0$). To see this we recall the Baston-Mason
result \cite{BasMas} stating that there are two neccessary conditions for a 
4-dimensional metric $g=g_{\mu\nu}\theta^\mu\theta^\nu$ to be
conformal to an Einstein metric. These, when expressed in terms of
the Cartan normal conformal connection $\tilde{\om}$, are\footnote{It
  is worthwile to note that for algebraically general metrics the
  Baston-Mason conditions $(i)-(ii)$ are also sufficient for the conformal
  Einstein property \cite{BasMas}.}
\be
(i)~~\der*\tilde{\Om}+\tilde{\om}\dz*\tilde{\Om}-*\tilde{\Om}\dz\tilde{\om}=0,~~~~{\rm and}~~~~(ii)~~[\tilde{\Om}^+_{~\mu\nu},\tilde{\Om}^-_{~\rho\sigma}]=0,\label{bama}
\ee
where
$\tilde{\Om}^\pm=\frac{1}{2}\tilde{\Om}^\pm_{~\mu\nu}\theta^\mu\dz\theta^\nu$.
Note that condition $(i)$ is equivalent to the vanishing of the Bach
tensor of $g$.\\

\noindent
Calculating $[\tilde{\Om}^-_{~32},\tilde{\Om}^+_{~31}]$ for the
Fefferman metrics (\ref{fef}) yields 
$$
[\tilde{\Om}^-_{~32},\tilde{\Om}^+_{~31}]=
i {\cal R}\bar{\cal R}\bma
0& 0& 0&1 & 0& 0\\
0& 0& 0& 0& 0& 0\\ 
0& 0& 0& 0& 0& 0\\ 
0& 0& 0&0& 0& 0\\
0& 0& 0& 0& 0& 1\\
0& 0& 0& 0& 0& 0
\ema,
$$
so that the above condition $(ii)$ is satisfied iff ${\cal
  R}=0$. This means that the corresponding
CR-structure is flat. It follows that if
${\cal R}=0$ the corresponding Fefferman metrics are
conformal to the Minkowski metric. In the non-flat (${\cal R}\neq
0$) case the Fefferman metrics are always not conformal to Einstein
metrics. Note also that despite of this fact the principal null
direction of the Fefferman metrics (which in the notation of (\ref{fef})
is tangent to the vector filed dual to the form $\theta^3$) is
geodesic and shear-free. It has nonvanishing twist and generates a
1-parameter conformal symmetry of $g_F$.\\

\noindent
We close this section with the formula for 
$\tilde{D}*\tilde{\Om}=\der*\tilde{\Om}+\tilde{\om}\dz*\tilde{\Om}-*\tilde{\Om}\dz\tilde{\om}$
which, for the Fefferman metrics (\ref{fef}), reads
\be
\tilde{D}*\tilde{\Om}=-\frac{2}{i}{\cal S}_1
\bma
0&0&0&1&0&0\\
0&0&0&0&0&0\\
0&0&0&0&0&0\\
0&0&0&0&0&0\\
0&0&0&0&0&1\\
0&0&0&0&0&0
\ema\theta^1\dz\theta^2\dz\theta^3,\label{23}
\ee
where ${\cal S}_1$ is defined by (\ref{ghj}).
This formula implies that the Fefferman metrics (\ref{fef}) satisfy
the Bach equations iff
\be
{\cal S}_1=0,\label{24}
\ee
or, what is the same,
$$
\der{\cal S}\dz\theta^2\dz\theta^3\dz(2\Om_2+3\bar{\Om}_2)\dz\bar{\Om}_3=0.
$$ 
The only known example of a CR-structure with a Fefferman metric
satisfying this condition is presented in \cite{NurPle}.

\section{Second order ODEs modulo point transformations}
\noindent
A second order ODE   
\be
\frac{\der^2 y}{\der x^2}~=~Q(~x,~y,~\frac{\der y}{\der x}~)\label{ode}
\ee
for a function $\R\ni x\to y=y(x)\in\R$, can be
alternatively written as a system of the two first order ODEs
$$
\frac{\der y}{\der x}~=~p,~~~~~~~~~~~~~~\frac{\der p}{\der x}~=~Q(x,y,p)
$$ 
for two functions $\R\ni x\to y=y(x)\in\R$ and  $\R\ni
x\to p=p(x)\in\R$. This system defines two (contact) 1-forms 
\be
\om^1=\der y-p\der x,~~~~~~~~~~~~~~\om^2=\der p-Q\der x,\label{forms}
\ee 
on a 3-dimensional manifold $J^1$, 
the {\it first jet space}, parametrized by
coordinates $(x,y,p)$. All the information about the ODE (\ref{ode})
is encoded in these two forms. For example, any solution  to
(\ref{ode}) is a curve $\gamma(x)=(~x,y(x),p(x)~)\subset J^1$ on
which the forms (\ref{forms}) vanish.\\

\noindent
The two contact 1-forms $(\om^1,\om^2)$ can be supplemented by
\be
\om^3=\der x,\label{form3}
\ee 
so that the three 1-forms $(\om^1,\om^2,\om^3)$ constitute a basis of
1-forms on $J^1$. This basis will be the basic object of
study in the following.\\

\noindent
Under the point transformation of
variables
$$
y\to \tilde{y}=\tilde{y}(x,y),~~~~~~~~~x\to \tilde{x}=\tilde{x}(x,y),
$$
the function $Q=Q(x,y,y')$ defining the differential equation
(\ref{ode}) changes in a rather
complicated way. The corresponding change of the basis
$(\om^1,\om^2,\om^3)$ is 
\be
\om^1\to\tilde{\om}^1=a_1\om^1,~~~~~~\om^2\to\tilde{\om}^2=a_2\om^2+a_3\om^1,~~~~~~\om^3\to\tilde{\om}^3=a_4\om^3+a_5\om^1,\label{tranreal}
\ee  
where $a_1,a_2,a_3,a_4,a_5$ are real functions on $J^1$ such that
$a_1a_2a_4\neq 0$ on $J^1$.\\

\noindent
It is now convenient, to introduce the following (a bit unusual)
notation. The reason for this will eventually become apparent.\\

\noindent
Let $i\neq 0$ denote a {\it real} number. In addition, let the {\it
  real} 1-forms $(\la,\mu,\bar{\mu})$ be defined by  
\be
\la=-i\om^1,~~~~~~~\mu=\om^2,~~~~~~~\bar{\mu}=\om^3.\label{lamubarmu}
\ee  
It follows from the definition of $(\om^1,\om^2,\om^3)$ that 
\be
\la\dz\mu\dz\bar{\mu}\neq 0,\label{f1}
\ee
\be
\der\la\dz\la\neq 0,\label{f2}
\ee
and that the forms $(\la,\mu,\bar{\mu})$ are given up to
transformations 
\be\la\to a\la,~~~~~~\mu\to b\mu+c\la,~~~~~~
\bar{\mu}\to\bar{b}\bar{\mu}+\bar{c}\la,\label{fri}\ee
with real functions $a,b,\bar{b},c,\bar{c}$ such that $ab\bar{b}\neq 0$.

\noindent
Conversely, given a 3-dimensional manifold $N$ equipped with three real
1-forms $(\la,\mu,\bar{\mu})$ satisfying
(\ref{f1})-(\ref{f2}) and defined up to transformations (\ref{fri}), 
we can associate with them a point equivalent class of a 2nd order ODE 
as follows. Since ${\rm dim} J^1=3$ we have
$$\der\la\dz\la\dz\bar{\mu}=0~~~~~{\rm and}~~~~~
\der\bar{\mu}\dz\la\dz\bar{\mu}=0.
$$
Hence the Fr\"{o}benius theorem \cite{Olver} applied to the forms
$\lambda,\bar{\mu}$ implies that there exist coordinates $(x,y,z)$ on $N$ such
that $\lambda=A\der x+B\der y$ and $\bar{\mu}=C\der x+H\der y$, where
$A,B,C,H$ are appropriate functions on $N$. Thus, 
modulo the freedom (\ref{fri}), the forms $\la$, $\bar{\mu}$ 
can be transformed to $\la=\der y-p\der x$, $\bar{\mu}=\der x$, where
$p$ is a certain function of coordinates $(x,y,z)$ on $N$. But, 
$0\neq \der\la\dz\la=\der p\dz\der y\dz\der x$ so $(x,y,p)$ can be
considered a new coordinate system on $N$. In this coordinates the form
$\mu$ can be written as $\mu=U\der x+V\der y+Z\der p$ so, by means of
transformations (\ref{fri}), can be reduced to $\mu=\der p- Q\der
x$ with $Q=Q(x,y,p)$ a certain real function on $N$. Thus, the original 
forms $(\la,\mu,\bar{\mu})$ define a point equivalent class of a
second order ODE $y''=Q(x,y,y')$. The above considerations prove the 
one-to-one correspondence between second order ODEs given modulo point
transformations and equivalence classes of the triples of real 1-forms
$(\la,\mu,\bar{\mu})$ on 3-manifolds satisfying (\ref{f1}),(\ref{f2})
and given up to (\ref{fri}). This enables us to reformulate an
equivalence problem for second order ODEs modulo point transformations
in much the same way as an equivalence problem for {\it
  nondegenerate} 3-dimensional CR-structures. 

\bd
Two second order ODEs, represented, by the
respective real 1-forms $(\la,\mu,\bar{\mu})$ and $(\la',\mu',\bar{\mu}')$, 
on the respective 3-manifolds $N$ and $N'$, 
are locally 
{\rm point equivalent}, iff there exist a local diffeomorphism 
$$
\phi:N\to N'
$$
and {\rm real} functions
$a\neq 0,b\neq 0,\bar{b}\neq 0,c,\bar{c}$ on $N$ such that
$$
\phi^*(\la')=a\la,~~~~~~~\phi^*(\mu')=b\mu+c\la,~~~~~~~
\phi^*(\bar{\mu})=\bar{b}\mu+\bar{c}\la.
$$
\ed

\noindent
This definition, when compared with Definition 2, indicates that we can
treat the forms $(\la,\mu,\bar{\mu})$ representing second order ODEs 
as the respective analogs of the forms 
$(\la,\mu,\bar{\mu})$ representing nondegenerate 3-dimensional
CR-structures. It also indicates that the solution for the equivalence 
problem for 2nd
order ODEs modulo point transformations should be given by a theorem
analogous to Theorem 1. Actually, with the above introduced notation,
in which all the three 1-forms
$(\la,\mu,\bar{\mu})$ are {\it real}, $i\neq 0$ is a {\it real} number and 
the `bar' symbol merely denotes that a given variable 
(a function, or a form) is totally {\it independent} of its non-bared 
counterpart, we obtain the solution of the equivalence problem for ODEs
by the following reinterpretation of Theorem 1. First, given a point 
equivalence class of 2nd order ODEs, represented by forms 
$(\la,\mu,\bar{\mu})$, we associate  with it the forms
\be
\theta^1=b\mu+c\la,~~~~\theta^2=\bar{b}\bar{\mu}+\bar{c}\la,~~~~
\theta^3=a\la.\label{osi}
\ee
Then the analog of Theorem 1 is as follows.  
\bt
Every second order ODE given modulo point transformations 
uniquely defines an 8-dimensional manifold $P$, {\rm real} 1-forms
$\theta^1,\theta^2,\theta^3,\Om_2,\bar{\Om}_2,\Om_3,\bar{\Om}_3,\Om_4$
and {\rm real} functions ${\cal R},\bar{\cal R},{\cal S},\bar{\cal S}$
on $P$ such that
\begin{itemize}
\item $\theta^1,\theta^2,\theta^3$ are as in (\ref{osi}),
\item $\theta^1\dz\theta^2\dz\theta^3\dz\Om_2\dz\bar{\Om}_2\dz\Om_3\dz\bar{\Om}_3\dz\Om_4\neq
  0$ at each point of $P$.
\end{itemize}
\noindent
The forms satisfy the following equations
\beq
&\der\theta^1=\Om_2\dz\theta^1+\Om_3\dz\theta^3\nonumber\\
&\der\theta^2=\bar{\Om}_2\dz\theta^2+\bar{\Om}_3\dz\theta^3\nonumber\\
&\der\theta^3=i\theta^1\dz\theta^2+(\Om_2+\bar{\Om}_2)\dz\theta^3\nonumber\\
&\der\Om_2=2
i\theta^1\dz\bar{\Om}_3+i\theta^2\dz\Om_3+\Om_4\dz\theta^3\nonumber\\
&\der\bar{\Om}_2=-2
i\theta^2\dz\Om_3-i\theta^1\dz\bar{\Om}_3+\Om_4\dz\theta^3\label{sys1}\\
&\der\Om_3=\Om_4\dz\theta^1+\Om_3\dz\bar{\Om}_2+{\cal
  R}\theta^2\dz\theta^3\nonumber\\
&\der\bar{\Om}_3=\Om_4\dz\theta^2+\bar{\Om}_3\dz\Om_2+\bar{\cal
  R}\theta^1\dz\theta^3\nonumber\\
&\der\Om_4=i\Om_3\dz\bar{\Om}_3+\Om_4\dz(\Om_2+\bar{\Om}_2)+\bar{\cal
  S}\theta^1\dz\theta^3+{\cal S}\theta^2\dz\theta^3\nonumber.
\eeq
The functions ${\cal R}$, $\bar{\cal R}$, ${\cal S}$, $\bar{\cal S}$, satisfy in addition 
\beq
&\der{\cal R}=-{\cal R}(\Om_2+3\bar{\Om}_2)-{\cal S}\theta^1+{\cal R}_2\theta^2+{\cal
  R}_3\theta^3\nonumber\\
&\der\bar{\cal R}=-\bar{\cal R}(\bar{\Om}_2+3\Om_2)-\bar{\cal S}\theta^2+\bar{\cal R}_2\theta^1+\bar{\cal
  R}_3\theta^3\label{funk1}
\eeq
and
\beq
&\der{\cal S}=-{\cal S}(2\Om_2+3\bar{\Om}_2)-i{\cal
  R}\bar{\Om}_3+{\cal S}_1\theta^1+{\cal S}_2\theta^2+{\cal S}_3\theta^3\nonumber\\
&\der\bar{\cal S}=-\bar{\cal S}(2\bar{\Om}_2+3\Om_2)+i\bar{\cal
  R}\Om_3+\bar{\cal S}_2\theta^1+{\cal S}_1\theta^2+\bar{\cal S}_3\theta^3,\nonumber
\eeq
with appropriate functions ${\cal R}_2,\bar{\cal R}_2,{\cal
  R}_3,\bar{\cal R}_3,{\cal S}_1,\bar{\cal S}_1,{\cal
  S}_2,\bar{\cal S}_2,{\cal S}_3,\bar{\cal S}_3$.\\
The function ${\cal S}_1$ satisfies 
\be 
{\cal S}_1=\bar{\cal S}_1.\label{bach1}
\ee
\et
Given an equation $y''=Q(x,y,y')$ and the standard coordinate system $(x,y,p)$
on $N=J^1$ we introduce the vector field 
$$
D=\partial_x+p\partial_y+Q\partial_p,~~~~~~~~~Q=Q(x,y,p).
$$  
The coordinates $(x,y,p)$ can be extended to a coordinate system
$(x,y,p,\rho,\phi,\gamma,\bar{\gamma},r)$ on $P$ 
in which the forms and functions of the above theorem can be written as follows:
\beq
&\la=-i(\der y-p\der x),~~~~~\mu=\der p-Q\der x,~~~~~\bar{\mu}=\der x\nonumber\\
&\nonumber\\
&\theta^1=\rho {\rm e}^{i\phi}(\mu+\ga \la),~~~~~
\theta^2=\rho {\rm e}^{-i\phi}(\bar{\mu}+\bar{\ga}
\la),~~~~~\theta^3=\rho^2\la\label{frc}\\
&\nonumber\\
&\Om_2=i\der\phi+\frac{\der\rho}{\rho}+\frac{1}{4i\rho^2}~[~6\ga\bar{\ga}i^2-6\bar{\ga}i
  Q_p-Q_{pp}-4ir\rho~]~\theta^3-\frac{2i\bar{\ga}}{\rho}{\rm
  e}^{-i\phi}\theta^1-\frac{{\rm e}^{i\phi}}{\rho}(i\ga-Q_p)\theta^2\nonumber\\
&\nonumber\\
&\bar{\Om}_2=-i\der\phi+\frac{\der\rho}{\rho}-\frac{1}{4i\rho^2}~[~6\ga\bar{\ga}i^2-2\bar{\ga}iQ_p-Q_{pp}+4ir\rho~]~\theta^3+\frac{i\bar{\ga}}{\rho}{\rm
  e}^{-i\phi}\theta^1+\frac{{\rm
    e}^{i\phi}}{\rho}(2i\ga-Q_p)\theta^2\nonumber\\
&\nonumber\\
&\Om_3=\frac{{\rm
    e}^{i\phi}}{\rho}~[~\der\ga-\frac{1}{6i^2\rho^2}(DQ_{pp}+6\ga^2\bar{\ga}i^3-6\ga\bar{\ga}i^2Q_p-3\ga i Q_{pp}-4Q_{py}-6\bar{\ga}iQ_y)\theta^3+\nonumber\\
&\frac{{\rm
      e}^{-i\phi}}{4i\rho}(2\ga\bar{\ga}i^2-2\bar{\ga}iQ_p-Q_{pp}-4ir\rho)\theta^1+\frac{{\rm e}^{i\phi}}{i\rho}(\ga^2i^2-\ga iQ_p-Q_y)\theta^2~]\nonumber\\
&\nonumber\\
&\bar{\Om}_3=\frac{{\rm
    e}^{-i\phi}}{\rho}~[~\der\bar{\ga}+\frac{1}{6i^2\rho^2}(6\ga\bar{\ga}^2 i^3-6\bar{\ga}^2i^2Q_p-3\bar{\ga} i Q_{pp}-Q_{ppp})\theta^3-\frac{{\rm
      e}^{-i\phi}}{\rho}\bar{\ga}^2i\theta^1-\nonumber\\
&\frac{{\rm e}^{i\phi}}{4i\rho}(2\ga\bar{\ga}i^2-2\bar{\ga}
  iQ_p-Q_{pp}+4ir\rho)\theta^2~]\nonumber\\
&\nonumber\\
&\Om_4=-\frac{i}{2\rho^2}\bar{\ga}\der\ga+\frac{1}{2\rho^2}(i\ga-Q_p)\der\bar{\ga}-\frac{\der
  r}{\rho}-\frac{r\der\rho}{\rho^2}+\nonumber\\
&\frac{1}{48i^2\rho^4}~[~8DQ_{ppp}-3Q_{pp}^2+8Q_pQ_{ppp}-12Q_{ppy}-12\ga
  iQ_{ppp}+\bar{\ga}(12iDQ_{pp}-24iQ_{py})-12\ga\bar{\ga}i^2Q_{pp}+\nonumber\\
&\bar{\ga}^2(24i^2DQ_p+12i^2Q_p^2-48i^2Q_y)-48i^3\ga\bar{\ga}^2Q_p+36\ga^2\bar{\ga}^2i^4+48\bar{\ga}r\rho
  i^2Q_p+48i^2\rho^2r^2~]~\theta^3-\nonumber\\
&\frac{{\rm
    e}^{i\phi}}{12i\rho^3}~[~6\ga\bar{\ga}^2i^3-6\bar{\ga}^2i^2Q_p+3\bar{\ga}iQ_{pp}-Q_{ppp}-12\bar{\ga}i^2r\rho~]~\theta^1-\nonumber\\
&\frac{{\rm e}^{i\phi}}{12i\rho^3}~[~DQ_{pp}-4Q_{py}-3i\ga
  Q_{pp}+6i\bar{\ga}(DQ_p-2Q_y)-6i^2\ga\bar{\ga}Q_p+6\ga^2\bar{\ga}i^3+12\ga i^2r\rho~]~\theta^2,\nonumber\\
&\nonumber\\
&{\cal R}=-\frac{{\rm e}^{2i\phi}}{6i^2\rho^4}w_1,~~~~~~~~~{\cal S}=-\frac{{\rm
    e}^{i\phi}}{3i^2\rho^5}~[\partial_pw_1+i\bar{\ga}w_1~],\nonumber\\
&\label{rr}\\
&\bar{\cal R}=-\frac{{\rm
    e}^{-2i\phi}}{6i^2\rho^4}w_2,~~~~~~~~~~~~~\bar{\cal S}=-\frac{{\rm e}^{-i\phi}}{3i^2\rho^5}~[~Dw_2+(2Q_p-i\ga)w_2~],\label{rrc}
\eeq
where we have introduced functions
$$
w_1=D^2Q_{pp}-4DQ_{py}-DQ_{pp}Q_p+4Q_pQ_{py}-3Q_{pp}Q_y+6Q_{yy},
$$
and
$$
w_2=Q_{pppp}
$$
which are the relative point invariants of the ODE.\\

\noindent
\noindent
Similarly as in the CR case, Theorem 2 can be reinterpreted in terms
of the language of Cartan connections. It follows that the manifold
$P$ of Theorem 2 is a Cartan bundle $H\to P\to J^1$, with $H$
a 5-dimensional parabolic subgroup of ${\bf SL}(3,{\bf R})$ group. The forms
$\theta^1,\theta^2,\theta^3,\Om_2,\bar{\Om}_2,\Om_3,\bar{\Om}_3,\Om_4$
of the theorem can be collected into a matrix of 1-forms 
$$
\omega=
\begin{pmatrix}
\frac{1}{3}(2\Om_2+\bar{\Om}_2)&i\bar{\Om}_3&-\frac{1}{2}\Om_4\\
&&\\
\theta^1&\frac{1}{3}(\bar{\Om}_2-\Om_2)&-\frac{1}{2}\Om_3\\
&&\\
2\theta^3&2i\theta^2&-\frac{1}{3}(2\bar{\Om}_2+\Om_2),
\end{pmatrix}
$$
which is now an ${\bf sl}(3,{\bf R})$-valued Cartan connection on $P$
(all the variables are real!). It follows from
equations (\ref{sys1}) that the curvature of
this connections is
$$
\Om=\der \om+\om\dz\om=
\begin{pmatrix}
0&i\bar{\cal R}\theta^1\dz\theta^3&-\frac{1}{2}\bar{\cal
  S}\theta^1\dz\theta^3-\frac{1}{2}{\cal S}\theta^2\dz\theta^3\\
&&\\
0&0&-\frac{1}{2}{\cal R}\theta^2\dz\theta^3\\
&&\\
0&0&0
\end{pmatrix}.
$$
It yields all the invariant information about the corresponding
point equivalent class of second order ODEs. In particular, the
ODEs corresponding to flat (${\cal R}=0$, $\bar{\cal R}=0$, ${\cal
  S}=0$, $\bar{\cal S}=0$) connections are given by the
conditions
$$w_1=0,~~~~~~~~w_2=0.$$
They are all point equivalent to the flat equation $y''=0$. We 
remark that the vanishing of $\cal R$ implies
vanishing of $\cal S$. Each of these two conditions is a point invariant
property of the corresponding ODE. However, the other pair of point
invariant conditions $\bar{\cal R}=0$, $\bar{\cal S}=0$ is totally
independent. This is the significant difference 
between the behaviour of CR structures and 2nd order ODEs. Indeed, in
the classification of nondegenerate 3-dimensional CR-structures there
are only two major branches: either ${\cal R}=0$ (in which case the
CR-structure is locally equivalent to ${\bf S}^3\subset {\bf C}^2$) or
${\cal R}\neq 0$. In the ODE case $\cal R$ and $\bar{\cal
  R}$ are unrelated and we have four main branches, corresponding to (i)
${\cal R}=0$, $\bar{\cal R}=0$, (ii) ${\cal R}=0$, $\bar{\cal R}\neq
0$, (ii') ${\cal R}\neq 0$, $\bar{\cal R}=0$, and (iii) ${\cal R}\neq
0$, $\bar{\cal R}\neq 0$. It follows that branches (ii) and (ii') are,
in a sense, dual to each other. To explain this duality we need to
introduce the Fefferman metric associated with an ODE.\\

\noindent
The system (\ref{sys1}) defining the invariant forms of Theorem
2 has all the qualitative properties of system (\ref{sys}) of Theorem
1. Thus, introducing the basis  
$X_1,X_2,X_3,Y_2,\bar{Y}_2,Y_3,\bar{Y}_3,Y_4$ of vector fields, 
the respective dual of forms $\theta^1,\theta^2,\theta^3,$
$\Om_2,\bar{\Om}_2,\Om_3,\bar{\Om}_3,\Om_4$, we see that the
distribution spanned by the four vector fields 
$Z_I=(Y_2+\bar{Y}_2,Y_3,\bar{Y}_3,Y_4)$ is integrable.
Moreover, the bilinear form 
$$
G=2\theta^1\theta^2+\frac{2}{3i}\theta^3(\Om_2-\bar{\Om}_2),
$$
which now has signature $(++--0000)$, has all $Z_I$s as degenerate
directions. This, when compared with the fact that $G$ is preserved up
to a scale during the Lie transport along $Z_I$s, shows that the
4-dimensional 
space $P/\hspace{-0.15cm}\sim$ of leaves of the distribution 
spanned by $Z_I$s is naturally equipped with the conformal class of
split signature metrics $[g_F]$, the
class to which the bilinear form $G$ naturally descends. We call the metrics
\be
g_F=2\theta^1\theta^2+\frac{2}{3i}\theta^3(\Om_2-\bar{\Om}_2)\label{fef1}
\ee
on $P/\hspace{-0.15cm}\sim$ the Fefferman metrics 
associated with a point equivalence class of a second order
ODE $y''=Q(x,y,y)$.\\

\noindent
The metrics $g_F$, when written in coordinates 
$(x,y,p,\rho,\phi,\gamma,\bar{\gamma},r)$ on $P$, read
\be
g_F=2\rho^2~[~(\der p-Q\der x)\der x-(\der y-p\der x)(\frac{2}{3}i\der
  \phi+\frac{2}{3}Q_p\der x+\frac{1}{6}Q_{pp}(\der y-p\der x))~].\label{ff1}
\ee
This enables us to (locally) identify the space 
parametrized by
$(x,y,p,\phi)$ with $P/\hspace{-0.18cm}\sim$ and the space parametrized by
$(x,y,p,\phi, \rho)$ with the space of all Fefferman metrics 
associated with a given $y''=Q(x,y,y')$.\\

\bp
The Fefferman conformal class of metrics $[g_F]$ associated with a
point equivalent class of ODEs has the following properties.
\begin{itemize}
\item Each $g_F$ has signature $(++--)$
\item The Weyl tensor of each $g_F$ has both, the self-dual and the 
anti-self-dual parts of Petrov type N. The self-dual part $C^+$ is 
proportional to ${\cal R}$ and the anti-self-dual part $C^-$ is
proportional to $\bar{\cal R}$.
\item $g_F$ satisfies the Baston-Mason conditions (\ref{bama}) 
if and only if the corresponding point equivalent class of equations
satisfies either ${\cal R}=0$ or $\bar{\cal R}=0$.
\end{itemize} 
\ep
The first two statements of the above proposoition are abvious in view
 of formulae (\ref{18})and (\ref{20}). To prove the last statement we 
calculate the Baston-Mason conditions (\ref{bama}) in coordinates 
$(x,y,p,\phi)$. A short calculation and the identity
\be
w_{1pp}=(D^2 + 3Q_pD + 2DQ_p + 2Q_p^2 -
Q_y)w_2\label{to1}
\ee
show that these conditions are equivalent to
\be
(i')~~~w_{1pp}=0~~~~~~~~~~~~~~~{\rm and}~~~~~~~~~~~~~~~(ii')~~~w_1 w_2=0,\label{to2}
\ee
where $(i')$ corresponds to the vanishing of the Yang-Mills current of 
the Cartan normal conformal connection $\tilde{\om}$ associated with
$g_F$ via (\ref{18}) and $(ii')$ corresponds to the Baston-Mason
condition $(ii)$ of (\ref{bama}). Comparing (\ref{to1}) and
(\ref{to2}) proves that the necessary and sufficient conditions for 
(\ref{bama}) are $w_1=0$ or
$w_2=0$. This in particular means that such metrics must be either
self-dual or anti-selfdual.\\

\noindent
{\bf Remark}\\
Note that one of the principal null directions of $g_F$,
generated by the vector field dual to the form $\theta^3$ of
(\ref{fef1}), is a conformal Killing vector field for $g_F$. 
It generates a congruence of null shear-free and twisting geodesics on
$P/\hspace{-0.15cm}\sim$. This statement, together with the above
Proposition 1 totally characterizes the Fefferman metrics $g_F$ \cite{Spar}.\\

\noindent
Since metrics $g_F$ are algebraically special 
(of type $N\times N'$ or its specializations) the Baston-Mason 
conditions (\ref{bama}) are not sufficient to guarantee the 
conformal Einstein property for them. All the Fefferman metrics which
are conformal to Einstein metrics are given in the Appendix.\\

\noindent
It is very easy to determine all classes of second order ODEs
correponding to the self-dual Fefferman metrics. These are all equations
for which $w_2=0$, i.e. all the equations of the form 
$$
y''=A_0(x,y)+A_1(x,y)y'+A_2(x,y)(y')^2+A_3(x,y)(y')^3,
$$
for which the function $Q$ is
an arbitrary polynomial
of the third order in the variable $p$. Finding classes of equations
corresponding to the anti-self-dual metrics is more difficult but,
surprisingly, possible, due to another notion of duality: the duality
of second order ODEs.\\

\noindent
Given a second order ODE in the form 
\be
\frac{\der^2 y}{\der x^2}=Q(x,y,\frac{\der y}{\der x})\label{du1}
\ee 
consider its general
solution $y=y(x,X,Y)$, where $X,Y$ are constants of integration. In
the space ${\bf R}^2\times {\bf R}^2$ parametrized by $(x,y,X,Y)$ this 
solution can be considered a 3-dimensional hypersurface 
$$N=\{(x,y,X,Y)\in {\bf R}^2\times {\bf R}^2~|~G(x,y,X,Y)=y-y(x,X,Y)=0~\}.$$
Assuming that $G(x,y,X,Y)=0$ can be solved with respect to $Y$ one
gets a function $Y=Y(X,x,y)$. Treating $x$ and $y$ as constant
parameters, we can eliminate them by double differentiation of $Y$
with respect to $X$. This means that $Y=Y(X,x,y)$
can be considered a solution to a second order ODE  
\be
\frac{\der^2 Y}{\der X^2}=q(X,Y,\frac{\der Y}{\der X}).\label{du2}
\ee
In a passage from (\ref{du1}) to (\ref{du2}) we have chosen $X$ to be
an independent variable and $Y$ to be the dependent one. But another
choices are possible. In general, we could have chosen two independent
functions $\xi=\xi(X,Y)$ and $\zeta=\zeta(X,Y)$ and have treated $\xi$
and $\zeta$ as an independent and dependent variables,
respectively. Then, after double differentiation with respect to
$\xi$, which eliminates the parameters
$(x,y)$, we would see that the function $\zeta$ also satisfied a
second order ODE, which would be quite different then (\ref{du2}). It
is, however, obvious that this other
second order ODE would be in the same point equivalence class as
(\ref{du2}). Thus, given a point equivalence class of ODEs generated
by (\ref{du1}) there is a uniquely defined equivalence class of ODEs
(\ref{du2}) associated with it. The class (\ref{du2}) is called a
class of the {\it dual} equations to the equations from the class
(\ref{du1}). The following example shows the usefulness of this concept.\\

\noindent
{\bf Example}\\
The relative invariants $w_1$ and $w_2$ calculated for the second order ODE 
\be
y''=\frac{a}{y^3},\label{geor}
\ee
where $a$ is a real constant, read
$$
w_1=\frac{72 a}{y^5}~~~~~~~~~~~~~{\rm and}~~~~~~~~~~~~w_2=0.
$$
Therefore the Fefferman metrics 
\be
g_F=2\rho^2~[~(\der p-\frac{a}{y^3}\der x)\der x-\frac{2}{3}i(\der
  y-p\der x)\der\phi~]\label{pg}
\ee
associated with the point equivalent class of ODEs generated by 
(\ref{geor}) are self-dual but not conformally
flat \footnote{These are actually the Sparling-Tod metrics well known in
  the twistor theory \cite{SparTod}}.\\  
The general solution $y=y(x,X,Y)$ of (\ref{geor}) depends on two
arbitrary constants of integration $(X,Y)$ and satisfies 
\be
y^2 = Y(x- X)^2+ \frac{a}{Y}.\label{sol11}
\ee
This generates a hypersurface  
$$
N=\{(x,y,X,Y)\in {\bf R}^2\times {\bf R}^2~|~y^2 = \frac{a}{Y} + 
Y(x-X)^2~\}
$$ 
in ${\bf R}^2\times {\bf R}^2$. Now, we will treat equation
(\ref{sol11}) as an equation for a function $Y=Y(X,x,y)$ of 
an independent variable $X$ and parametrized by $x$ and $y$. 
Differentiating (\ref{sol11}) with respect to $X$ and keeping $x$ and
$y$ constants we get
$$0=Y'(x-X)^2-2Y(x-X)-a\frac{Y'}{Y^2}.$$
Solving for $x$ and differentiating once more with respect to $X$ we
find that $Y=Y(X)$
satisfies the second order ODE
$$
Y''=-\frac{- Y^4 Y'^2 +a Y'^4-2 Y^2 Y'^2\sqrt{Y^4+a Y'^2}}
{Y^5+Y^3\sqrt{Y^4+a Y'^2}}.
$$
This is an equation which generates the point equivalence class of
equations dual to (\ref{geor}). This equation has $w_2\neq 0$ for each
$a\neq 0$. A direct (but lenghty!) calculation shows that 
$$
q=q(X,Y,P)=-\frac{- Y^4 P^2 +a P^4-2 Y^2 P^2\sqrt{Y^4+a P^2}}
{Y^5+Y^3\sqrt{Y^4+a P^2}}
$$
has $w_1=0$. This is a general fact known already to Elie Cartan \cite{Cartpc}. More
formally, we have the following Proposition.
\bp 
The point equivalence class of dual ODEs to a point equivalence class 
of 2nd order ODEs for which $w_2=0$ and $w_1\neq 0$ has $w_2\neq 0$
and $w_1=0$. 
\ep
We have already noted (in the Example above) 
that this proposition enables one to find nontrivial 
solutions to a quite complicated differential equation $w_1=0$. Note
also that applying the Proposition one can obtain quite nontrivial
anti-self-dual metrics from a rather dull ones (Calculate
the Fefferman metrics for $q$ of the Example, and compare it with
(\ref{pg})). Finally, note that it follows from the Proposition that
in the classification scheme of the 2nd order ODEs modulo
point transformations, the classification of the $w_1=0$ and $w_2\neq 0$
case can be obtained from the classification of the simpler $w_2=0$
and $w_1\neq 0$ case.\\ 

\noindent
In Ref. \cite{Cartpc} Proposition 2
is only briefly mentioned\footnote{For an unexperienced reader it is rather
  hard to find a trace of the Proposition in the text of Ref. \cite{Cartpc}. We are very
  grateful to Mike Crampin \cite{Crampin} for clarifying this point
  for us. We also take this opportunity to present his understanding of
  the last paragraph of Cartan's paper \cite{Cartpc}. Due to
  equations (\ref{frc}), (\ref{rr})-(\ref{rrc}), which represent the 
transformation
  properties of the forms $(\om^1,\om^2,\om^3)$ and the relative 
invariants $w_1$ and $w_2$, if both $w_1$ and $w_2$ are nonvanishing, 
the following forms
  are (modulo sign) well defined on $J^1$ 
$$I_1=(w_1w_2)^{\frac{1}{4}}\om^1,~~~~~~~~I_2=(w_1w_2)^{\frac{1}{2}}\om^1\dz\om^2\dz\om^3$$
$$I_3=w_1^{\frac{1}{8}}w_2^{\frac{5}{8}}
\om^1\dz\om^2,~~~~~~~~I_4=w_1^{\frac{5}{8}}w_2^{\frac{1}{8}}\om^1\dz\om^3.$$}. It could be proven by the following line of argument.\\

\noindent
The switch between the dual equations $y''=Q(x,y,y')$ and
$Y''=q(X,Y,Y')$ essentially means the switch between the contact forms
$\om^2=\der p-Q\der x$ and $\om^3=\der x$. To see this consider the
general solution $y=y(x,X,Y)$ of the original equation. Now, pass from
the canonical coordinates $(x,y,p)$ on the first jet bundle $J^1$ to
the new coordinates $(s,X,Y)$ defined via 
$$x=s,~~~~~~~y=y(s,X,Y)~~~~~~{\rm and}~~~~~~p=y_s(s,X,Y).$$ 
In these new coordinates the contact forms
$$\om^1=\der y-p\der x,~~~~~~\om^2=\der p-Q\der x,~~~~~~~\om^3=\der
x$$
associated with the original equation are given by 
$$\om^1=y_X\der X+y_Y\der Y,~~~~~~\om^2=y_{sX}\der X+y_{sY}\der
Y,~~~~~~~\om^3=\der s.$$
Thus in the point equivalence class of the forms $\om^1$ and $\om^2$ there are
forms 
$$\om^1=\der Y+\frac{y_X}{y_Y}\der X~~~~~{\rm and}~~~~\om^2=\der X.$$
Moreover the
condition $0\neq\der\om^1\dz\om^1=\der(\frac{y_X}{y_Y})\dz\der X\dz\der
Y$ implies that the three functions 
$$X,~~~~~Y~~~~~{\rm and}~~~~~P=-\frac{y_X}{y_Y}$$
constitute a coordinate system on $J^1$. Therefore $s$ must be a
function of these three variables: $s=s(X,Y,P)$. This means that in the
equivalence class of the form 
$$\om^3=\der s=s_X\der X+s_Y\der Y+s_P\der P$$
there is $\om^3$ which can be writen as 
$$\om^3=\der P-q(X,Y,P)\der X,~~~~~~~~~~~~~{\rm where}~~~~~~~~~~~~
q=-\frac{s_X+s_YP}{s_{P}}.$$
Summarizing, we are able to
introduce a coordinate system $(X,Y,P)$ on the first jet bundle $J^1$
in which the point equivalence class of the contact forms associated
with the original equation $y''=Q(x,y,y')$ can be written as
$$\om^1=\der Y-P\der X,~~~~~~~\om^2=\der X,~~~~~~~\om^3=\der
P-q(X,Y,P)\der X.$$
But this enables us to interpret the $(X,Y,P)$ coordinates as cannonical
coordinates for the first jet bundle associated with the contact
forms 
$$\om^1=\der Y-P\der X,~~~~~~~\om^3=\der P-q(X,Y,P)\der X,~~~~~~~\om^2=\der X$$
of the differential equation $Y''=q(X,Y,Y')$. Because of the original
definitions of $X$ and $Y$ this is clearly the dual equation to 
$y''=Q(x,y,y')$.  Note that this interpretation requires the switch
between the forms $\om^2$ and $\om^3$. Note also that this switch is 
compatible with transformations (\ref{tranreal}) which treat $\om^2$
and $\om^3$ on the equal footing mixing each of them with $\om^1$ only. \\

\noindent
Once the switch between $y''=Q(x,y,y')$ and $Y''=q(X,Y,Y')$ is
understood from the point of view of the switch between $\om^2$ and 
$\om^3$ it is easy to see that the passage from a differential
equation to its dual changes the
role of the invariants $w_1$ and $w_2$. Indeed, looking at the
curvature $\Om$ of the Cartan connection $\om$ associated with the
equation $y''=Q(x,y,y')$ we see that the invariant $w_1$ is associated 
with the $\om^3\dz\om^1$ term 
and the invariant $w_2$ is associated with the $\om^2\dz\om^1$
term. Thus, the switch between $\om^2$ and $\om^3$ caused by the
switch between the equation and its dual, amounts in the switch of 
$w_1$ and $w_2$.\\

\noindent
From the geometric point of view the switch between the mutually dual
2nd order ODEs can be understood as a transformation that interchanges 
two naturally defined congruences of lines on the jet bundle
$J^1$. This bundle is naturally fibred over $J^0$ - the plane 
parametrized by $(x,y)$. The fibres of $J^1\to J^0$ are 1-dimensional 
and, in the natural coordinates $(x,y,p)$ on $J^1$, can be 
specified by fixing $x$ and $y$. They 
generate the first congruence of lines on $J^1$. The
other congruence is defined by the point equivalence class of the 
equation $y''=Q(x,y,y')$ in the following way. The equation
$y''=Q(x,y,y')$ equippes $J^1$ with
the total differential vector field $X_{\bar{\mu}}=D$.  Any other equation from the
point equivalence class of $y''=Q(x,y,y')$ defines the
total differential that differs from $D$ by a scaling functional factor.
Thus the lines tangent to all of these total diferentials are well
defined on $J^1$ and generate the second natural congruence of lines.
Each of the above congruences on $J^1$ defines a natural direction
of vector fields tangent to them but, in the canonical coordinates 
$(x,y,p)$ on $J^1$, only the lines of the first congruence can be
defined as lines tangent to a particularly simple vector field 
$X_\mu=\partial_p$. The passage from the equation to its dual changes
the picture: it switches between $X_\mu$ and $X_{\bar{\mu}}$, so that
the jet bundle $J^1$ is now interpretd as a bundle with 1-dimensional
fibres tangent to $X_{\bar{\mu}}=\partial_P$. The space of such fibres may then
be identified with the plane parametrized by $(X,Y)$ and the
congruence tangent to $X_\mu$ by the congruence tangent to the total
differential of the dual equation.

\section{Realification of a 3-dimensional CR-structure}
The analogy between 3-dimensional CR-structures and 2nd order ODEs
described in the previous two sections can be used to associate a
point equivalent class of 2nd order ODEs with a 3-dimensional
CR-structure. This may introduce a new insight in General Relativity,
where the 3-dimensional CR structures correspond \cite{tafel,trautopt,traut}
to congruences of
shear-free and null geodesics in space-times\footnote{Such space-times, the Lorentzian 4-dimensional
manifolds admiting a null congruence of shear-free geodesics, are
called {\it Robinson manifolds} \cite{NurTraut,t1,t2} and are known to be
the analogs of Hermitian manifolds of four dimensional Riemannian geometry.}.
In particular, the
shear-free congruence of null geodesics associated with the celebrated
Kerr space-time, can be interpreted in terms of a certain class of
point equivalent second order ODEs. In this section we provide a
framework for this kind of considerations concentrating on
3-dimensional CR-structures that admit 2-dimensional group of local
symmetries\footnote{Note that the Kerr congruence is in this class of
  examples \cite{NPracmag}}.\\

\noindent
A 3-dimensional CR-structure with two symmetries can be locally
described by real coordinates $(u,x,y)$ in which a representative
of the forms from the class $[(\la,\mu)])$ is given by
\be
\la=\der u+f(y)\der x,~~~~~~\mu=\der x+i\der y,~~~~~~\bar{\mu}=\der x-i\der y.\label{kerr}
\ee
Here the function $f=f(y)$ is real and $i^2=-1$.\\

\noindent
To pass from the above CR-structures to the corresponding point
equivalence class of 2nd order ODEs we require that the symbol $i$ is a 
{\it  nonvanishing real constant} so that the forms (\ref{kerr}) become {\it
  real} and we can interpret them as the forms that via (\ref{fri})
define a point equivalence class of 2nd order ODEs. To find a
particular representative of on ODE in this class we introduce new
coordinates $(\bar{u},\bar{x},\bar{y})$, which are related to
$(u,x,y)$ by
$$
u=\bar{u},~~~~~~y=\bar{y},~~~~~~x=\bar{x}+i\bar{y}.
$$   
Since $i$ is now {\it real} and nonvanishing this is a {\it real} 
transformation of the coordinates. It brings $(\la,\mu,\bar{\mu})$
to the form
$$
\la=\der (~\bar{u}+i\hspace{-0.15cm}\int\hspace{-0.15cm}f(\bar{y}~)
\der\bar{y}~)~+~f(\bar{y})\der\bar{x},~~~~~~~~~~~~~~
\mu=\der\bar{x}+2i\der\bar{y},~~~~~~~~~~~~~~~\bar{\mu}=\der\bar{x}.
$$ 
After another coordinate transformation
$$X=\bar{x},~~~~~Y=\bar{u}+i\hspace{-0.15cm}\int\hspace{-0.15cm}f(\bar{y})\der\bar{y},~~~~~~P=-f(\bar{y})$$
and an application of the chain rule one sees
that in the class (\ref{fri}) of 1-forms $(\la,\mu,\bar{\mu})$ there
are forms
$$\la=\der Y-P\der X,~~~~~~\mu=\der
P~-~\frac{1}{2i}~~f'(y)_|\hspace{-0.11cm}{_{_{_{|y=f^{-1}(-P)}}}}\hspace{-0.1cm}\der X,~~~~~~\bar{\mu}=\der X.$$
This means thet the point equivalence classes of 2nd order ODEs
associated with the CR-structures generated by $f=f(y)$ are represented by the
equations of the form
$$Y''=\frac{1}{2i}f'(y)_|\hspace{-0.11cm}{_{_{_{|y=f^{-1}(-Y')}}}}.$$
Consider, in particular, the family of CR-structures which have 
3-dimensional symmetries of Bianchi type $VI_k$ \cite{NT}. They are
parametrized by the Bianchi type parameter $k<0$ and are represented
by the function $f=y^n$, where $n=n(k)$ is an apropriate \cite{NT}
function of $k$. Then the point equivalent class of ODEs associated with each $n$
is given by
\be
Y''=\frac{n}{2i}(-Y')^{1-\frac{1}{n}}.\label{rcr}
\ee
The application of the above results to the new understanding
of the congruences of shear-free and null geodesics in space-time and,
in particular, to the Kerr congruence\footnote{This congruence may be 
represented by $f=\frac{1}{{\rm
  cosh}^2(y)}$ \cite{NPracmag}.} may be of some use in General Relativity Theory. Here we
focus on properties of Fefferman metrics associated with equations
(\ref{rcr}). It is known \cite{NurPle} that the Fefferman
metrics associated with the 3-dimensional CR-structures admitting
three symmetries of Bianchi type corresponding to $n=-3$ satisfy the
Bach equations. They are the only known {\it Lorentzian} metrics of twisting
type $N$ which satisfy the Bach equations and which are not conformal to
Einstein metrics. Similarly the Fefferman metrics for the equations
point equivalent to the equation (\ref{rcr}) with $n=-3$ provide an example of
{\it split} signature metrics of type $N\times N'$ which satisfy the Bach
equations and which are never conformal to Einstein
metrics. Explicitely, these metrics are conformal to 
$$g_F=2~(~i P^{\frac{2}{3}}\der P+\frac{3}{2}P^2\der X~)~\der X-2~(~\der
Y-P\der X~)~(~\frac{2}{3}i^2 P^{\frac{2}{3}}\der\Phi-\frac{1}{9}\der
  Y-\frac{11}{9}P\der X~).$$
These metrics can be generalized by considering 2nd order ODEs point
equivalent to 
$$
y''=h(y'),
$$
where $h=h(p)$ is a sufficiently smooth real function. Then the
Fefferman metrics associated with such equations are given by
\be
g_F=2\rho^2~[~(\der p-h\der x)\der x-(\der y-p\der x)(\frac{2}{3}i\der
  \phi+\frac{2}{3}h'\der x+\frac{1}{6}h''(\der y-p\der x))~].\label{spf}
\ee
The invariants $w_1$ and $w_2$ for these metrics are 
$$
w_1=h^2 h'''',~~~~~~~~~~~~w_2=h''''
$$
so the metrics are not conformal to Einstein metrics iff 
$$
h''''\neq 0.
$$
This, together with (\ref{to2}), shows that every solution to the equation 
$$
h^2 h''''=ap+b,
$$
where $a$ and $b$ are real constants such that at least one of them
does not vanish defines, via
(\ref{spf}), a {\it split} signature conformal class of metrics that
satisfy the Bach equations, are of type $N\times N'$ and are not
conformal to Einstein metrics. By an appropriate 
complexfication of these metrics one may get the generalization of 
the Lorentzian Bach non-Einstein metrics (45) of Ref. \cite{NurPle}.
\section{Acknowledgements}
This work is a byproduct of the lectures on `Cartan's equivalence
method' which one of us (PN) delivered at the Department of Physics
and Astronomy of
Pittsburgh Unversity in the winter of 2002-3. The topic presented here
would never
even have been touched by us if Ezra Ted Newman would not requested to 
hear such lectures. He, his collaborators and the PhD students of
the Mathematics and Physics $\&$ Astronomy Departments of Pittsburgh 
University were our first audience and our inspiration.\\

\noindent
We are very grateful to Lionel Mason, David Robinson and Andrzej Trautman, who
carefully read the draft of these paper suggesting several
crucial improvements.\\

\noindent   
The paper was completed during the stay of one of us at King's College
London. We acknowledge support from EPSRC grant GR/534304/01 while at
King's College London, NSF
grant PHY-0088951 while at Pittsburgh University and the KBN grant 2
P03B 12724 while at Warsaw University.

\section{Appendix}
In this Appendix we find all Fefferman metrics (\ref{ff1}) which are 
conformally equivalent to the Einstein
metrics. First, we write them in the form
$$
g_F=2\theta^1\theta^2+\theta^3\theta^4,
$$
with
\beq
&\theta^1=\rho (\der p-Q\der x)\nonumber\\
&\theta^2=\rho \der x\label{ff2}\\
&\theta^3=-\rho^2 (\der y-p\der x)\nonumber\\
&\theta^4=\frac{2}{3}i\der\phi+\frac{2}{3}Q_p\der
x+\frac{1}{6}Q_{pp}(\der y-p\der x),\nonumber\\
&~\nonumber
\eeq
where $\rho=\rho(x,y,p,\phi)$ is a function on $J^1$. We search for all
$\rho$ and $Q=Q(x,y,p)$ for which the Einstein
equations 
\be
Ric (g_F)=\Lambda g_F\label{einsla}
\ee
\noindent
are satisfied. Assuming that $g_F$ obey equations (\ref{einsla}) we
use the implied Weyl tensor identity 
$
\nabla_\mu C^\mu_{~\nu\rho\sigma}=0.
$
Its $\{113\}$ and $\{223\}$ components imply the equations 
\be
(3\rho_\phi-i \rho)w_2=0~~~~~~{\rm and}~~~~~~~(3\rho_\phi+i \rho)w_1=0.
\ee
Thus, the metrics $g_F$ may be conformal to Einstein metrics only in
the following three cases:
\begin{itemize}
\item[($\bullet$)] $w_1=w_2=0$
\item[($\bullet\bullet$)] $w_2=0$, $w_1\neq 0$ and $\rho=A(x,y,p){\rm
  exp}(-\frac{i\phi}{3})$
\item[($\bullet\bullet\bullet$)] $w_2\neq 0$, $w_1= 0$ and $\rho=A(x,y,p){\rm exp}(\frac{i\phi}{3})$
\end{itemize}
\vskip 1.cm
\noindent
Case ($\bullet$) corresponds to ODEs with flat Cartan
connection. These have conformally flat Fefferman metrics.\\ 

\noindent
In the ($\bullet\bullet$) case
the general solution for the Einstein equations (\ref{einsla}) is
given by (\ref{ff2}), where either 
\be
\rho=\frac{{\rm exp}(-\frac{2}{3}a)}{p+2b}{\rm exp}(-\frac{i\phi}{3}),
\ee
and
\be
Q=p^3 c+p^2 (6bc-2a_y)+p(12b^2c-3b_y-6ba_y-a_x)+8b^3c-2bb_y-4b^2a_y-2b_x-2ba_x,
\ee
or 
\be
\rho={\rm exp}(-\frac{2}{3}a){\rm exp}(-\frac{i\phi}{3}),
\ee
and
\be
Q=p^2 a_y+2pa_x+b,
\ee
with $a=a(x,y)$, $b=b(x,y)$, $c=c(x,y)$ arbitrary functions of
variables $x,y$. All these solutions are Ricci flat. They exhaust the
list of Fefferman metrics which have $w_2=0$ and are conformal to non-flat
Einstein metrics. Since $Q$ in these solutions depend
only on at most {\it three} arbitrary functions of two variables and generic 
$Q$ for which $w_2=0$ depends on {\it four} functions of two
variables, then not all ODEs with $w_2=0$ have Fefferman metrics 
which are conformal to Einstein metrics.\\

\noindent
The Einstein equations in the ($\bullet\bullet\bullet$) case reduce to 
\be
\rho=a {\rm exp}(\frac{i\phi}{3}),\label{ff3}
\ee
where the function $a=a(x,y)$ satsifies a single differential equation
\be
36(Da)^2+6a(-3D^2a+Q_pDa)+a^2(6DQ_p-18Q_y-4Q_p^2)=0.\label{rownanko}
\ee
It follows, as it
should be, that one of the
integrability conditions for this equation is $w_1=0$. It is not clear
what other conditions should be imposed on $Q$ to guarantee
the solvability of (\ref{rownanko}). If, for example, $Q$ satisfies  
$$
6DQ_p-18Q_y-4Q_p^2=0,
$$
a simple solution is given by $a=const$.\\

\noindent
Finally, we remark that all
the metrics $g_F$ which, via (\ref{ff2}), correspond to the solutions
(\ref{ff3})-(\ref{rownanko}) are Ricci flat. This, when compared
with the results of case ($\bullet\bullet$), proves that
conformally non flat Fefferman-Einstein metrics have vanishing
cosmological constant.\\

\end{document}